\documentclass[a4paper,12pt,amsmath]{article}

\usepackage{theorem, euscript}
\usepackage{amsmath, amsfonts, amssymb, graphics, mathrsfs}
\usepackage{enumerate}
\font\tenbb=msbm10 \font\sevenbb=msbm7 \font\fivebb=msbm5
\newfam\bbfam
\textfont\bbfam=\tenbb \scriptfont\bbfam=\sevenbb
\scriptscriptfont\bbfam=\fivebb

\theoremstyle{remark}
\newtheorem{theorem}{Theorem}

\newtheorem{proposition}[theorem]{Proposition}
\newtheorem{lemma}[theorem]{Lemma}

\newtheorem{remark}[theorem]{Remark}

\author{ }
\newenvironment{Proof}[1][D{\'e}monstration]{\begin{trivlist}
\item[\hskip \labelsep {\bfseries #1}]}{\flushright
$\Box$\end{trivlist}}

\newcommand{\qed}{\nobreak \ifvmode \relax \else
      \ifdim\lastskip<1.5em \hskip-\lastskip
      \hskip1.5em plus0em minus0.5em \fi \nobreak
      \vrule height0.75em width0.5em depth0.25em\fi}

\renewcommand{\det}{{\text{d{\'e}t}}}

\def\Ac{{\cal A}}

\def\Ac{{\cal A}}

\def\g{\gamma}

\def\part{\partial}

\def\text{\hbox}

\def\build#1_#2^#3{\mathrel{
\mathop{\kern 0pt#1}\limits_{#2}^{#3}}}

\font\tenbb=msbm10 \font\sevenbb=msbm7 \font\fivebb=msbm5
\newfam\bbfam
\textfont\bbfam=\tenbb \scriptfont\bbfam=\sevenbb
\scriptscriptfont\bbfam=\fivebb

\def\Ac{{\cal A}}

\def\g{\gamma}

\def\build#1_#2^#3{\mathrel{
\mathop{\kern 0pt#1}\limits_{#2}^{#3}}}

\def\Im{\mathop{\rm Im}\nolimits}

\def\det{\mathop{\rm det}\nolimits}

\def\log{\mathop{\rm log}\nolimits}

\def\calh{{\mathcal H}}

\def\g{\gamma}

\def\part{\partial}

\def\text{\hbox}

\def\zm{\begin{Proof}}

\def\endzm{\end{Proof}}

\makeindex

\title{Rankin-Cohen Deformations and Representation Theory}
\author{Yi-Jun YAO}

\begin{document}
\maketitle
\begin{abstract}In this paper, we use the unitary
representation theory of $SL_2(\mathbb R)$ to understand the
Rankin-Cohen brackets for modular forms.  Then we use this
interpretation to study the corresponding deformation problems
that Paula Cohen, Yuri Manin and Don Zagier initiated. Two
uniqueness results are established.
\end{abstract}

\section{Introduction}

Let $\Gamma$ be a congruence subgroup of $PSL(2, {\mathbb Z})$.
For $k\in{\mathbb N}$, a modular form of weight $2k$ is a complex
function $f$ on the upper half plane $\mathbb H$ which
satisfies(\cite{Houches89}):

\begin{itemize}
\item  (holomorphy) $f$ is holomorphic.
\item (modularitity) For $\gamma=\displaystyle \left(\begin{array}{cc}a & b\\
                                                                    c &
d\end{array}\right)\in
                                                                    \Gamma$
                                                                    and $z\in\mathbb
                                                                    H$,
                                                                    $f\Big|_{2k}\gamma=f$,
                                                                    where
\begin{equation}
\left(f\Big|_{2k}\gamma\right)(z)=(cz+d)^{-2k}f\left(\frac{az+b}{cz+d}\right),
\end{equation}

\item (growth condition at the boundary) We ask that $|f(z)|$ would be controlled by a polynomial in $\max\{1, \Im(z)^{-1}\}$.
\end{itemize}

We note by ${\mathcal M}(\Gamma)=\displaystyle
\bigoplus_{k\in\mathbb{N}} {\mathcal M}_{2k}(\Gamma)$ the graded
algebra (by the weight) of modular forms with respect to this
group.

In the 50's Rankin began the study of bidifferential
operators over ${\mathcal M}(\Gamma)$ which produce
new modular forms,  and twenty years later
Henri Cohen gave a complete answer (cf. \cite{C75}) by proving that
all these operators are linear combinations of the following brackets

\begin{equation}
[f, g]_n=\sum_{r=0}^n (-1)^r {n+2k-1 \choose n-r}{n+2l-1 \choose
r} f^{(r)}g^{(n-r)} \in {\mathcal M}_{2k+2l+2n}(\Gamma),
\end{equation}

\noindent where $f\in {\mathcal M}_{2k}$ and $g\in {\mathcal
M}_{2l}$ are two modular forms, and
$f^{(r)}=\displaystyle\left(\frac{1}{2\pi
i}\frac{\partial}{\partial z}\right)^r f$.

These brackets attracted interest of several authors. In \cite{Z94},  Zagier used the Ramanujan derivation $X:{\mathcal
M}_{2k} \rightarrow {\mathcal M}_{2k+2}$:

\begin{equation}
X f= \frac{1}{2  \pi i} \,  \frac{df}{dz} - \frac{1}{2  \pi i} \,
\frac{\partial}{\partial z} (\log \eta^4) \cdot k f.
\end{equation}

\noindent and introduced two series of elements by induction:

\begin{equation}
f_{r+1} =\partial f_r + r(r+2k-1) \Phi f_{r-1} \, \,  ,  \, \,
g_{s+1} =\partial g_s + s(s+2l-1) \Phi g_{s-1},
\end{equation}

\noindent where $\Phi=\displaystyle\frac{1}{144} E_4\in\mathcal{M}_4$ and $E_4$ is the Eisenstein series of weight 4. He showed that

\begin{equation}\label{CanRC}
\sum_{r=0}^n (-1)^r {n+2k-1 \choose n-r}{n+2l-1 \choose r}
f_{r}g_{n-r} = [f, g]_n,
\end{equation}
which turned the modularity of $[f, g]_n$ obvious as all the $f_r$
and $g_{n-r}$ are modular.

Moreover, he showed that for all associative $\mathbb Z$ (or
$\mathbb N$)-graded  algebra having a derivation which increase
the degree by 2, and for all element $\Phi$ of degree 4, the
formula (\ref{CanRC}) defines a {\it canonical Rankin-Cohen
algebra} structure.

\begin{remark}
When $\Phi=0$, the situation is simplified to what Zagier called
{\sl standard  Rankin-Cohen algebra}.
\end{remark}

\begin{remark}
We remind the readers that in the above definitions only the
modularity is used, so we can do the same for nonholomorphic
functions.
\end{remark}

About the same time, Paula Cohen,  Yuri Manin and Don Zagier
established a bijective correspondance between the modular forms
and the invariant formal pseudodifferential operators, They showed
that the following formula (plus linear extension) defines an
associative product over ${\mathcal M}(\Gamma)[[\hbar]]$: for two
modular forms $f\in{\mathcal M}_{2k}$,  $g\in{\mathcal M}_{2l}$,

\begin{equation}
\mu^\kappa (f, g):=\sum_{n=0}^\infty t_n^\kappa(k, l) [f, g]_n,
\end{equation}

\noindent where the coefficients are given by

\begin{equation}\label{mukappa}
t_n^\kappa(k, l)=\left(-\frac{1}{4}\right)^n \sum_{j\geq 0}
{n\choose 2j}\frac{ \displaystyle{-\frac{1}{2}\choose
j}{\kappa-\frac{3}{2}\choose j}{\frac{1}{2}-\kappa\choose j}}
{\displaystyle{-k-\frac{1}{2}\choose j}{-l-\frac{1}{2}\choose
j}{n+k+l-\frac{3}{2}\choose j}}.
\end{equation}

A special case is when $\kappa=\frac{1}{2}$ or $\frac{3}{2}$, and the product is reduced to what
Eholzer claimed to be an associative product

\begin{equation}\label{ehol}
f\star g :=\sum_{n=0}^\infty [f, g]_n.
\end{equation}

\begin{remark}
In this formulation, only the modularity of $f$ is used, we do not
need neither holomorphy, nor the growth condition near the
boundary.
\end{remark}

In 2003, Connes and Moscovici related the Hopf algebra $\calh_1$
introduced their study of transversal index theory, which governs
the local symmetry in calculating the index of a transversal
elliptic operator, to the Rankin-Cohen brackets. By taking into
account the work of Cohen-Manin-Zagier and of Eholzer, especially
(\ref{ehol}), Connes and  Moscovici proved a theorem stating that
for every $\calh_1$ on an algebra $\Ac$ with certain extra
structure, there exists a familly of formal deformations of $\Ac$
where the general terms of the deformed products are defined by
some generalized Rankin-Cohen brackets (\cite{CM03-2})

In a joint work with P.Bieliavsky and X.Tang(\cite{BTY}), we have
studied the deformation question from a quite different point of
view. We used the deformation quantization theory of Fedosov to
construct a realization of Rankin-Cohen deformations. More
precisely we found a specific symplectic connection on the upper
half plane and on the corresponding Weyl algebra we found the same
induction relation as that of Connes-Moscovici while calculating
the deformed product. Then by an analoguous argument, we re-obtain
the above theorem of Connes-Moscovici.

\bigskip
In this paper, we study the brackets via the unitary representation theory of $SL_2({\mathbb R})$
and then apply the results thus obtained to the deformation questions.

The rest of this paper is organized as follows: first a
(relatively) explicit interpretation of the Rankin-Cohen brackets
are given via the representation theory of $SL_2({\mathbb R})$.
The principal result is the following theorem: \footnote{The
specialists in the domain certainly know this long time ago, as
showed by a remark that Deligne made in 1973 (cf. Remark
\ref{deligne}),  but before finish writing my Ph.D. thesis(Nov.
2006), I had not found any detailed  presentation of the this
result.}
\smallskip

\noindent {\bf Theorem}. {\it Let $f\in{\mathcal M}_{2k},
g\in{\mathcal M}_{2l}$ be two modular forms. Let
 $\pi_f\cong \pi_{\deg f}, \pi_g\cong \pi_{\deg
g}$ be the corresponding discrete series representations of
$SL_2({\mathbb R})$. The tensor product of these two
representations can be decomposed into a direct sum of discrete
series,

\begin{equation}
\pi_f\otimes \pi_g =\bigoplus_{n=0} \pi_{\deg f+\deg g+2n}.
\end{equation}

\noindent The Rankin-Cohen bracket $[f, g]_n$ gives (up to scale)
the vectors of minimal $K$-weight in the representation space of
the component $\pi_{\deg f+\deg g+2n}$;}

These representations are constructed in the following way: let
$f\in {\mathcal M}_{2k}(\Gamma)$ be a modular form, we associate
to it a function on $\Gamma\backslash SL_2({\mathbb R})$ by using
the following map: for $g=\left(\begin{array}{cr}
                           a & b\\
                            c & d
                            \end{array}\right)\in SL_2({\mathbb
                            R})$,

\begin{equation}
(\sigma_{2k} f)(g)=f\Big|_{k}g (i)=(ci+d)^{-2k}
f\left(\frac{ai+b}{ci+d}\right).
\end{equation}

This function belongs to

$$C^\infty(\Gamma\backslash SL_2({\mathbb R}), 2k)=\{F\in
C^\infty(\Gamma\backslash SL_2({\mathbb R})), F(gr_\theta)= \exp(
i 2k\theta) F(g) \}.$$ By taking into account the natural right
action of $SL_2({\mathbb R})$ on $C^\infty(\Gamma\backslash
SL_2({\mathbb R}))$:

\begin{equation}
(\pi(h) F)(g)=F(gh),
\end{equation}

\noindent we obtain a representation of $SL_2({\mathbb R})$ and so
of the complexified Lie algebra ${\mathfrak sl}_2({\mathbb C})$ by
taking the smallest invariant subspace which contains the orbit of
$\sigma_{2k} f$. We show that this representation is a discrete
series of weight $2k$. In the end, we pull all the vectors in a
basis of the representation space back to a subspace of $C^\infty(
{\mathbb H})$ by using the inverse of the $\sigma_{2(k+n)}$'s,
$n\geq 0$.

\smallskip

Then we use this representation theory interpretation to study
certain properties of the deformed products, and mainly we can get
the next two results:

\smallskip

\noindent{\bf Theorem}. {\it Cohen-Manin-Zagier have found {\sl
all} formal deformed associative products
$\ast:\widetilde{\mathcal M}[[\hbar]]\times\widetilde{\mathcal
M}[[\hbar]]\rightarrow \widetilde{\mathcal M}[[\hbar]]$ defined by
linear extension and the formula

\begin{eqnarray}
f\ast g & = & \sum \frac{A_n(\deg f,  \deg g)}{(\deg f)_n (\deg
g)_n} [f, g]_n \hbar^n,
\end{eqnarray}

\noindent where $\widetilde{\mathcal M}$ is the space of functions
which satisfy the
 modularity condition,  and the notation
$(\alpha)_n:=\alpha(\alpha+1)\cdots(\alpha+n-1)$.  We ask morevoer
$A_0=1$ and $A_1(x, y)=xy$. }

\noindent{\bf Proposition}. {\it Let $\Gamma$ be a congruence
subgroup of $SL_2({\mathbb Z})$ such that ${\mathcal M}(\Gamma)$
admits the unique factorization property (for example
$SL_2({\mathbb Z})$ itself), let $F_1$, $F_2$, $G_1$, $G_2\in
{\mathcal M}(\Gamma)$ such that

\begin{equation}
RC(F_1, G_1)=RC(F_2, G_2),
\end{equation}

\noindent as formal series in ${\mathcal M}(\Gamma)[[\hbar]]$,
then there exists a constant $C$ such that}

\begin{equation}
F_1=C F_2,  G_2= C G_1.
\end{equation}
{\bf Acknowledgement}: The author would like to thank his Ph.D.
advisor Alain Connes for guiding him into this interesting area and
for his constant support. He wants to thank Don Zagier for his
inspiring course given at Coll\`ege de France and to Henri
Moscovici, Fr\'ed\'eric Paugam, Jean-Pierre Labesse and Fran\c{c}ois
Martin for having helped him on several important points.

\section{ From modular forms to discrete series}

In this part we will describe another way to understand these
Rankin-Cohen brackets. We will partially follow the argument that
Jean-Pierre Labesse indicated (\cite{Lab05}):

\medskip

Let $f\in {\mathcal M}_{2k} (\Gamma)$ be a modular form of weight $2k$
with respect to a congruence subgroup $\Gamma$ of $SL_2(\mathbb Z)$. We will associate
 a  $\Gamma$-invariant function over $\Gamma\backslash
SL_2({\mathbb R})$ to it.

\medskip

We define

\begin{equation}
(\sigma_{2k} f)(g)=f\Big|_{k} g (i)=(ci+d)^{-2k}
f\left(\frac{ai+b}{ci+d}\right), \, \, \, \,
\end{equation}

\noindent for $g=\left(\begin{array}{cr}
                           a & b\\
                            c & d
                            \end{array}\right)\in SL_2({\mathbb R})$.
This function is invariant under the left
translation of the group $\Gamma$: let
$\g\in \Gamma$,  $f|_k \g g= (f|_k \g)|_k g=f|_k g$.

We verify also that for

\begin{equation}
r_\theta=\left(\begin{array}{cr}
                           \cos\theta & \sin\theta\\
                            -\sin\theta & \cos\theta
                            \end{array}\right)\in SL_2({\mathbb
                            R}),
\end{equation}

\noindent we have

\begin{eqnarray}\label{eqn1.0}
(\sigma_{2k} f)(gr_\theta) &=& \exp(i 2k\theta) (\sigma_{2k}
f)(g).
\end{eqnarray}

\medskip

 In fact, $\sigma_{2k}$ gives a bijection between

\begin{equation}
C^\infty(\Gamma\backslash {\mathbb H},  2k)=\left\{F\in
C^\infty({\mathbb H}), f(\gamma.z)= (cz+d)^{2k}f(z),
\gamma=\left(\begin{array}{cr}
                           a & b\\
                            c & d
                            \end{array}\right)\in \Gamma \right\}.
\end{equation}

\noindent and

\begin{equation}
C^\infty(\Gamma\backslash SL_2({\mathbb R}), 2k)=\{F\in
C^\infty(\Gamma\backslash SL_2({\mathbb R})), F(gr_\theta)= \exp(
i 2k\theta) F(g) \}.
\end{equation}

Take the space of smooth functions $C^\infty(\Gamma\backslash
SL_2({\mathbb R}))$,  we have a natural right action of
$SL_2({\mathbb R})$ on $\Gamma\backslash SL_2({\mathbb R})$: for
$F\in C^\infty(\Gamma\backslash SL_2({\mathbb R}))$,

\begin{equation}
(\pi(h) F)(g)=F(gh).
\end{equation}

We take the smallest invariant subspace under the action of
$SL_2({\mathbb R})$ which contains the orbit of $\sigma_{2k} f$ for a form
 $f\in {\mathcal M}_{2k}$,  and we are interested in the action of Lie algebra
${\mathfrak sl}_2({\mathbb R})$ on this space. We adopt the notations that S. Lang use in his book \cite{Lang75}.
A basis of this Lie algebra is

\begin{equation}
V=\left(\begin{array}{cc}
           0 & 1 \\
            1 & 0 \end{array}\right), \, \,
            H=\left(\begin{array}{cc}
           1 & 0 \\
            0 & -1 \end{array}\right), \, \,
            W=\left(\begin{array}{cc}
           0 & 1 \\
            -1 & 0 \end{array}\right),
\end{equation}

\noindent while a basis for the complexified Lie algebra ${\mathfrak
sl}_2({\mathbb C})$ is

\begin{equation}
E_+=\left(\begin{array}{cc}
           1 & i \\
            i & -1 \end{array}\right), \, \,
            E_-=\left(\begin{array}{cc}
           1 & -i \\
            -i & -1 \end{array}\right), \, \,
            W=\left(\begin{array}{cc}
           0 & 1 \\
            -1 & 0 \end{array}\right),
\end{equation}

\begin{eqnarray}
\exp(tV)&=&\left(\begin{array}{cc}
           \cosh t & \sinh t \\
            \sinh t & \cosh t \end{array}\right),
            \exp(tH)=\left(\begin{array}{cc}
           \exp t & 0 \\
            0 & \exp(-t) \end{array}\right), \cr
\exp(tE_+)&=&\left(\begin{array}{cc}
           1+t & it \\
            it & 1-t \end{array}\right),
            \exp(tE_-)=\left(\begin{array}{cc}
           1+t & -it \\
            -it & 1-t \end{array}\right), \cr
            \exp(tW)&=&\left(\begin{array}{cc}
           \cos t & \sin t \\
            -\sin t & \cos t \end{array}\right).
\end{eqnarray}

Now we take an arbitrary holomorphic function $\xi$ over the upper
half plane $\mathbb H$, for {\it all $k$}, we define

\[(F_k \xi) (g):= (\sigma_{2k} \xi)(g).\]

We calculate first the action of the base vectors described above on $F_k \xi$. We find

\begin{eqnarray}
(L_V F_k \xi) (g)  &=&
(-2k)\frac{di+c}{ci+d} (F_k \xi)(g )+2\left(F_{k+1}
\frac{d\xi}{dz}\right)(g), \cr (L_H F_k \xi)(g)  & = & (-2k)\frac{ci-d}{ci+d} (F_k
\xi)(g)+2i \left(F_{k+1} \frac{d\xi}{dz}\right)(g),
\end{eqnarray}

\noindent which implies

\begin{eqnarray}
L_{E_+}(F_{k} \xi)(g) & = &
2\left[(-2k)\frac{ci-d}{ci+d} (F_{k} \xi)(g) + 2i \left(F_{k+1}
\frac{d\xi}{dz}\right)(g)\right], \cr L_{E_-}(F_{k} \xi)(g)& = &
(L_H-iL_V)(F_{k} \xi)(g)=0.
\end{eqnarray}

\noindent And we have also

\begin{eqnarray}
 (L_W F_k \xi) (g) & = & 2ki (\sigma_{2k} \xi)(g) = 2ki (F_k \xi)(g).
\end{eqnarray}

  So by induction, we have

\begin{lemma}  For $n\in{\mathbb N}$,
\begin{enumerate}
\item $ \displaystyle (L_{E_+})^n (F_{k} \xi) = 2^n \sum_{t=0}^n
(-1)^{n-t} \frac{n!}{t!}{2k+n-1\choose n-t}
\left(\frac{ci-d}{ci+d} \right)^{n-t}(2i)^t\left(F_{k+t} \frac{d^t
\xi}{dz^t}\right)(g)$;
\item $L_W (L_{E_+})^n (F_{k} \xi)(g)=(2k+2n)i (L_{E_+})^n (F_{k}
\xi)(g)$;
\item $L_{E_-} (L_{E_+})^n (F_{k} \xi)(g)=-4n(2k+n-1) (L_{E_+})^{n-1} (F_{k}
\xi)(g)$.
\end{enumerate}
\end{lemma}

Next we calculate the action of the Casimir operator defined by

\begin{equation}\label{Casdef}
\omega = V^2 + H^2 - W^2=\frac{1}{2}(E_+ E_- + E_- E_+)-W^2.
\end{equation}

\noindent The above calculation shows that for each vector
$(L_{E_+})^n F_k \xi$

\begin{eqnarray}\label{Cascal}
\omega (L_{E_+})^n F_k \xi & = &
\frac{1}{2}[-4n(2k+n-1)-4(n+1)(2k+n)](L_{E_+})^n F_k \xi
+(2k+2n)^2(L_{E_+})^n F_k \xi \cr &= &4k(k-1)(L_{E_+})^n F_k \xi.
\end{eqnarray}

Thus the Casimir acts on the space generated by the $(L_{E_+})^n
F_k \xi$'s as constant.
\smallskip

\noindent If we start by a modular form $f$ (so a holomorphic
function) of weight $2k$ and form a vector space generated by the
functions $(L_{E_+})^n F_k f$. The above argument shows then
${\mathfrak sl}_2({\mathbb C})$ also acts and the Casimir acts as
the multiplication by the constant $4k^2-4k$. So we have a
representation of ${\mathfrak sl}_2({\mathbb C})$.

   Now we prove its irreducibility: for all operator  $T$ which commutes with the representation,  $[T,
E_-]=0$ implies that for the vector of minimal weight $F_k f$,  $T
F_k f$ is still a vector of minimal weight (for it's sent to zero by
 $E_-$),  so there is a constant $\lambda$ such that $T F_k f=\lambda F_k f$.
By the same argument, by $E_- T (E_+ F_k
f)= T  E_- (E_+ F_k f) = T(8k F_k f)=8k\lambda T_k f$,  we have
$T (E_+ F_k f)=\lambda E_+ F_k f$. So by induction we show that  $T$ acts by constant, the  representation is therefore irreducible. The representation theory of
$SL_2({\mathbb R})$ implies (cf. \cite{Lang75}, \cite{Vo}):

\begin{proposition}
What we have constructed is an irreducible representation of the Lie algebra
${\mathfrak sl}_2({\mathbb C})$ which is the infinitesimale version of the discrete series of the group
$SL_2({\mathbb R})$ of weight $2k$.
\end{proposition}
\medskip

When we take all these functions of  $C^\infty(SL_2({\mathbb R}))$
back to the space $C^\infty( {\mathbb H})$ by using the
bijectivity of the maps $\sigma_{2k+2n}$,  we get a representation
of ${\mathfrak sl}_2({\mathbb C})$, denoted by $\pi_f$. We denote
by $E_+, E_-, W$ the operators which correspond to $L_{E_+},
L_{E_-},  L_W$.

First by

\begin{eqnarray}
(\sigma_{2k+2} E_+ f) (g)&=&L_{E_+}(\sigma_{2k} f)(g)\cr & = &
2\left[(-2k)\frac{ci-d}{ci+d} (\sigma_{2k} f)(g) + 2i
\left(\sigma_{2k+2} \frac{df}{dz}\right)(g)\right]\cr &=&
2\left[2k\frac{1}{\Im \displaystyle\frac{ai+b}{ci+d}}\sigma_{2k+2}
f + 2i \sigma_{2k+2} \frac{df}{dz}\right](g),
\end{eqnarray}

\noindent we can define

\begin{equation}\label{deftildepar}
\widetilde{X} f :=-\frac{1}{8\pi}(E_+)f= \frac{1}{2\pi
i}\frac{df}{dz} - \frac{2kf}{4\pi Im z},
\end{equation}

\noindent which is called Shimura operator by some authors and
played an important role in Henri Cohen's paper \cite{C75} .

 In fact we can verify directly that

\begin{lemma}
Let $f$  be a differentiable function such that

\begin{equation*}
f\left(\frac{az+b}{cz+d}\right)=(cz+d)^{2k} f(z),
\end{equation*}

\noindent we have,

\begin{equation*}
\widetilde{X} f\left(\frac{az+b}{cz+d}\right)=(cz+d)^{2k+2}
\widetilde{X}f(z).
\end{equation*}
\end{lemma}

\noindent{\bf Proof.} It's sufficient to use

\begin{equation}
\Im \left(\frac{az+b}{cz+d}\right)=\Im
\left(\frac{az+b}{cz+d}\cdot \frac{c\bar{z}+d}{c\bar{z}+d}
\right)=\frac{\Im z}{|cz+d|^2}.
\end{equation}

The claim can be obtained by the following calculation:

\begin{eqnarray*}
\hspace{-0.5cm}\widetilde{X}
f\left(\frac{az+b}{cz+d}\right)&=&\frac{1}{2\pi
i}\frac{\partial}{\partial
z}\left(f\left(\frac{az+b}{cz+d}\right)\right)\Big\slash
\frac{\partial}{\partial
z}\left(\frac{az+b}{cz+d}\right)-\frac{2k}{4\pi \Im
\left(\frac{az+b}{cz+d}\right)}f\left(\frac{az+b}{cz+d}\right)\cr
\hspace{-0.5cm}&=& \frac{1}{2\pi
i}\left[(cz+d)^{2k}\frac{df}{dz}+2k
(cz+d)^{2k-1}f(z)\right](cz+d)^2\cr\hspace{-0.5cm} & &
\hspace{1.5cm} -\frac{2k}{4\pi}\frac{(cz+d)(c\bar{z}+d)}{\Im
z}(cz+d)^{2k} f(z)\cr \hspace{-0.5cm}&=&
(cz+d)^{2k+2}\frac{1}{2\pi
i}\frac{df}{dz}+(cz+d)^{2k+1}\frac{2k}{4\pi\Im
z}(cz+d)f(z)\cr\hspace{-0.5cm} &=& (cz+d)^{2k+2}
\widetilde{X}f(z). \Box
\end{eqnarray*}

By reiterating this operation, we get the following
correspondence:

\begin{eqnarray*}
  \left(-\frac{1}{8\pi}\right)^n\frac{1}{2k\cdots (2k+n-1)} (E_+)^n f & \leftrightarrow & \frac{1}{2k\cdots (2k+n-1)} \left( \frac{1}{2\pi i}\frac{\partial}{\partial z}- \frac{Y}{2\pi Im z}\right)^n f .
\end{eqnarray*}

\noindent where $Yf=kf$ is the Euler operator. Using the
representation theory of $SL_2({\mathbb R})$, we can choose the
vectors on the right hand side to form a basis, we can write

\begin{equation}
\varphi_n = \frac{1}{2k\cdots (2k+n-1)} \left( \frac{1}{2\pi
i}\frac{\partial}{\partial z}- \frac{Y}{2\pi Im z}\right)^n f,
\end{equation}

\noindent for $n\in{\mathbb N}$. The action of the Lie algebra
${\mathfrak sl}_2({\mathbb C})$ is given by

\begin{eqnarray}
E_+ \varphi_n &=& (-8\pi)(2k+n)\varphi_{n+1}, \\
E_- \varphi_n &=&\frac{n}{2\pi} \varphi_{n-1}, \\
W \varphi_n &=& 2ni \varphi_n.
\end{eqnarray}

We introduce an operator $\widetilde{\partial}$ such that
$\widetilde{\partial} \varphi_n=\varphi_{n+1}$, then

\begin{equation}\label{definitiontildepar}
\varphi_n=\widetilde{\partial}^n \varphi_0=\widetilde{\partial}^n
f.
\end{equation}

And moreover we have

\begin{lemma}
Let $f$ be a smooth function which satisfies the modularity
condition of weight $2k$, then,

\begin{equation}\label{der}
f^{(m)} :=\left(\frac{1}{2\pi i}\frac{\partial}{\partial
z}\right)^{m}f= m! \sum_{r=0}^m \frac{1}{(4\pi
y)^r}\frac{\widetilde{X} ^{m-r}}{(m-r)!} {2k+m-1 \choose r} f.
\end{equation}
 \end{lemma}

This implies exactly

\begin{equation}\label{combi}
[f, g]_n =\sum_{r=0}^n (-1)^r \widetilde{X}^r {2k+n-1\choose n-r}
f \widetilde{X}^{n-r} {2l+n-1\choose r}g,
\end{equation}

\noindent for $f\in {\mathcal M}_{2k}, g \in {\mathcal M}_{2l}$,
because

\begin{eqnarray*}
[f, g]_n & = & \sum_{r=0}^n (-1)^r {n+2k-1 \choose n-r}{n+2l-1
\choose r} f^{(r)}g^{(n-r)} \cr
 & = & \sum_{r=0}^n (-1)^r {n+2k-1 \choose n-r}{n+2l-1 \choose r}\cr
 & & \hspace{1cm} \left(r!\sum_{s=0}^r \frac{1}{(4\pi y)^s}{2k+r-1\choose s}\frac{\widetilde{X}^{r-s}}{(r-s)!}f\right)
 \cr & & \hspace{2cm}\left((n-r)!\sum_{t=0}^{n-r} \frac{1}{(4\pi y)^t}{2l+n-r-1\choose
t}\frac{\widetilde{X}^{n-r-t}}{(n-r-t)!}g\right)\cr & = & \sum_{s,
t} \frac{1}{(4\pi y)^{s+t}} \Bigg( \sum_{r=s}^{n-t} (-1)^r {n+2k-1
\choose n-r}{n+2l-1 \choose
r}\frac{r!}{(r-s)!}\frac{(n-r)!}{(n-r-t)!}\cr & &
\hspace{3cm}{2k+r-1\choose s}{2l+n-r-1\choose t}
\widetilde{X}^{r-s} f \widetilde{X}^{n-r-t} g \Bigg)
\end{eqnarray*}

It's clear that when $u=s+t,  v=r-s$ (and so $n-r-t=n-u-v$) are
all fixed,  the coefficient of $\widetilde{X}^{v} f
\widetilde{X}^{n-u-v} g$ is

\begin{eqnarray*}
& &\displaystyle\sum_s (-1)^{s+v} {n+2k-1 \choose n-v-s}{n+2l-1
\choose v+s}\cr & & \hspace{1cm}
\displaystyle\frac{(v+s)!}{v!}\frac{(n-v-s)!}{(n-v-u)!}{2k+v+s-1\choose
s}{2l+n-v-s-1\choose u-s}\cr & = & \displaystyle (-1)^{v} \sum_s
(-1)^s
\frac{(n+2k-1)!}{(2k+v+s-1)!(n-v-s)!}\frac{(n+2l-1)!}{(2l+n-v-s-1)!(v+s)!}
\cr & & \hspace{1cm}
\displaystyle\frac{(v+s)!}{v!}\frac{(n-v-s)!}{(n-v-u)!}\frac{(2k+v+s-1)!}{
s!(2k+v-1)!}\frac{(2l+n-v-s-1)!}{(u-s)!(2l+n-u-v-1)!}\cr & = &
\displaystyle (-1)^{v} \sum_s (-1)^s
\displaystyle\frac{(n+2k-1)!(n+2l-1)!}{(2k+v-1)!v!(n-v-u)!(2l+n-u-v-1)!}\frac{1}{s!(u-s)!}\cr
&=& (-1)^{v}
\frac{(n+2k-1)!(n+2l-1)!}{(2k+v-1)!v!(n-v-u)!(2l+n-u-v-1)!u!}\sum_s
(-1)^s \displaystyle\frac{u!}{s!(u-s)!},
\end{eqnarray*}

\noindent which is non-zero if and only if $u=0$, i.e., $s=t=0$.
We get thus the result. $\Box$

\medskip

We will see immediately a more conceptual explanation of this
identity.

\section{Construction of the brackets}

\medskip

Given two representations of $SL_2({\mathbb R})$(and the
corresponding derived representation of ${\mathfrak sl}_2({\mathbb
R})$ or ${\mathfrak sl}_2({\mathbb C})$), we're interested in
their tensor product. In fact, we have the following theorem of J.
Repka (cf. \cite{Rep78}):

\begin{theorem}
For two discrete series of $SL_2({\mathbb R})$, their tensor
product has the following decomposition:(for $m, n\geq 1$)

\begin{equation}
\pi_m\otimes \pi_n \cong \pi_{m+n}\oplus \pi_{m+n+2}\oplus
\pi_{m+n+4}\oplus \cdots \cong \bigoplus_{k=0}^\infty
\pi_{n+m+2k}.
\end{equation}
\end{theorem}

\medskip

To adapt this theorem (Lie algebra version) into our situation, we
give a special consideration on the representation space. More
precisely,

\begin{proposition}\label{propasso}
Given two modular forms $f\in{\mathcal M}_{2k}$, $g\in {\mathcal
M}_{2l}$, then in the decomposition

\begin{equation}
\pi_f \otimes \pi_g =\bigoplus_{n=0} \pi_{\deg f+ \deg g+2n},
\end{equation}

\noindent a vector of minimal $K$-weight of $\pi_{\deg f+ \deg
g+2n}$ has the form

\begin{eqnarray}\label{tilpartilx}
&&\frac{1}{n!}\sum_{r=0} (-1)^r {n \choose r}
\widetilde{\partial}^r f\otimes \widetilde{\partial}^{n-r} g\cr &
=&\frac{1}{(2k)_n(2l)_n}\sum_{r=0}^n (-1)^r \widetilde{X}^r
{2k+n-1\choose n-r} f \otimes \widetilde{X}^{n-r} {2l+n-1\choose
r}g.
\end{eqnarray}

Under the map defined by the product:

\begin{equation}\label{scmdef}
{\textsc m}: f\otimes g \longmapsto fg,
\end{equation}
this corresponds to a modular form of weight $2k+2l+2n$ which can
be expressed as

\begin{equation}
\frac{1}{2k(2k+1)\cdots (2k+n-1)2l(2l+1)\cdots (2l+n-1)}[f,
g]_n=\frac{1}{(2k)_n(2l)_n}[f, g]_n.
\end{equation}
\end{proposition}

\noindent \bf Proof. \rm The first part is a consequence of the
fact that the space of minimal $K$-weight vectors is exactly the
kernel of the operator $\Delta E_-=E_-\otimes 1+ 1\otimes E_-$, we
have

\begin{eqnarray*}
& & \Delta E_-\left(\sum_{r=0} (-1)^r {n \choose r}
\widetilde{\partial}^r f\otimes \widetilde{\partial}^{n-r}
g\right)\cr &=& \sum_{r=0} (-1)^r {n \choose r}
\left(E_-(\widetilde{\partial}^r f)\otimes
\widetilde{\partial}^{n-r} g+\widetilde{\partial}^r f\otimes
E_-(\widetilde{\partial}^{n-r} g)\right)\cr &=&
\frac{1}{2\pi}\sum_{r=0} (-1)^r {n \choose r}
\left(r\widetilde{\partial}^{r-1} f\otimes
\widetilde{\partial}^{n-r} g+(n-r)\widetilde{\partial}^r f\otimes
\widetilde{\partial}^{n-r-1} g\right)\cr &=&
\frac{1}{2\pi}\sum_{r=0} \left((-1)^r {n \choose
r}(n-r)+(-1)^{r+1}{n\choose r+1}(r+1)\right)
\widetilde{\partial}^r f\otimes \widetilde{\partial}^{n-r-1} g\cr
&=& 0.
\end{eqnarray*}

The second half is just (\ref{combi}). The operator $\textsc m$ is
a twister between the subrepresentation in the tensor product and
the representation constructed from  $[f, g]_n$. $\Box$

\medskip

\noindent {\bf N.B. In this construction, we can only determine
the coefficients up to scale.}

Furthermore, the formulation of Rankin-Cohen brackets using the
operator $\widetilde{X}$ can be naturally generalized to all pair
of functions $(f, g)\in \widetilde{\mathcal M}^2$, where

\begin{equation}
\widetilde{\mathcal M}(\Gamma) := \bigoplus_{k}
\widetilde{\mathcal M}_{2k} (\Gamma) := \bigoplus_{k} \left\{ f:
{\mathbb H}\rightarrow {\mathbb C},  f\Big|_{2k}\gamma=f, \, \,
\forall \gamma\in\Gamma \right\}
\end{equation}

\noindent is the space of smooth complex functions on the upper
half plane which satisfy (only) the modularity condition.

But in this case we do not have a general discrete series
interpretation as above.

\begin{remark}\label{deligne}

In fact, the relation between the tensor products of discrete
series representations and Rankin-Cohen brackets was already
observed some 35 years ago as one can find the following remark of
P.Deligne made in 1973 \cite{D73}: there he talked about discrete
series of $GL(2)$:

 ``Remarque 2.1.4. L'espace $F(G,  GL(2, {\mathbb Z}))$ ci-dessus est stable par produit. D'autre
part,  $D_{k-1}\otimes D_{l-1}$ contient les $D_{k+l+2m} (m \geq
0)$ . Pour $m = 0$,  ceci correspond au fait que le produit $fg$
d'une forme modulaire holomorphe de poids $k$ par une de poids
$l$,  en est une de poids $k+l$ . Pour $m = 1$,  en coordonn\'ees
(1.5.2)(remark: this should be 1.1.5.2),  on trouve que
$\displaystyle l \frac{\partial f}{\partial z}.g -
kf.\frac{\partial g}{\partial z}$ est modulaire holomorphe de
poids $k+l+2$,  et ainsi de suite. De m\^eme dans le cadre
ad\'elique.''

In fact, here what we get is the modularity of $\displaystyle
\frac{1}{k} \frac{\partial f}{\partial z}.g -
f.\frac{1}{l}\frac{\partial g}{\partial z}$.

After the main part of the paper was written (as one chapter of my
thesis in French), M. Weissman(\cite{Weissman06})posted on Arxiv a
paper which is along the line of Deligne's remark.

\end{remark}

\begin{remark}
We notice also there is an interpretation of these Rankin-Cohen
brackets Using the theory of transvectants. Especially in a recent
paper (\cite{Grad}),  El Gradechi treated the Rankin-Cohen
brackets in a very similar way as we did above.

\end{remark}

\section{Applications to Formal Deformations}

In this part, we study the formal deformations constructed from
the Rankin-Cohen brackets, more precisely we are interested in the
products $\ast: \widetilde{\mathcal M}(\Gamma) [[\hbar]]\times
\widetilde{\mathcal M}(\Gamma) [[\hbar]]\rightarrow
\widetilde{\mathcal M}(\Gamma) [[\hbar]]$ defined by linearity and
the formula:

\begin{eqnarray}\label{prodformfaib}
f\ast g & = & \sum \frac{A_n(\deg f, \deg g)}{(\deg f)_n(\deg
g)_n} \left(\sum_{r=0}^n (-1)^r \widetilde{X}^r {2k+n-1\choose
n-r} f \widetilde{X}^{n-r} {2l+n-1\choose r}g\right) \hbar^n\cr
&=& \sum \frac{A_n(\deg f,  \deg g)}{(\deg f)_n (\deg g)_n} [f,
g]_n \hbar^n,
\end{eqnarray}

\noindent where $f, g\in \widetilde{\mathcal M}$. We ask
furthermore $A_0=1$ et $A_1(x, y)=xy$. The main concern is to have
an associative product. First we have

\begin{proposition} \label{fort=faible}
If the $A_n$'s give rise to an associative product, then in the
expansion of $(f\ast g)\ast h$ and $f\ast (g\ast h)$, the
coefficients of every $\displaystyle \widetilde{X}^r f
\widetilde{X}^s g \widetilde{X}^t h$ are the same.
\end{proposition}

\noindent{\bf Proof.} In fact, we only need to show the equality
of the coefficients for $\displaystyle \frac{\partial^r
f}{\partial z^r}\frac{\partial^s g}{\partial z^s}\frac{\partial^t
h}{\partial z^t}$ and we prove this by contradiction. Assume that
there are functions $f_0, g_0, h_0 \in \widetilde{\mathcal M}$ and
an index triple $(r_0, s_0, t_0)$ such that the coefficient of
$\displaystyle \frac{\partial^{r_0} f}{\partial
z^{r_0}}\frac{\partial^{s_0} g}{\partial
z^{s_0}}\frac{\partial^{t_0} h}{\partial z^{t_0}}$ in $(f_0\ast
g_0)\ast h_0 -f_0\ast (g_0\ast h_0)$ is non-zero. So the
associativity of the product $\ast$ gives rise to a differential
equation which is satisfied by {\it all} $f\in{\mathcal M}_{\deg
f_0}, g\in{\mathcal M}_{\deg g_0}, h\in{\mathcal M}_{\deg h_0}$.

Now the only constrain on these functions are their invariance
under the action of $\Gamma$, which implies that we have the
freedom to modify the functions in the {\it interior} of a
fundamental domain. So in a small open set contained in the
fundamental domain, we can have some $f_1, g_1, h_1$ such that
$\displaystyle \frac{\partial^r f_1}{\partial
z^r}=\frac{\partial^s g_1}{\partial z^s}=\frac{\partial^t
h_1}{\partial z^t}=0$,  $0\leq r, s, t\leq n$, $r\neq r_0,  s\neq
s_0,  t\neq t_0$; and
\[\frac{\partial^{r_0} f}{\partial
z^{r_0}}\neq 0, \, \, \frac{\partial^{s_0} g}{\partial
z^{s_0}}\neq 0, \, \,  \frac{\partial^{t_0} h}{\partial
z^{t_0}}\neq 0.\]

But this gives us a contradiction. The proposition is
then proved. $\Box$

For three functions $f, g$ and $h$ in $\widetilde{\mathcal M}$,
the objects $(f\ast g)\ast h$ and $f\ast (g\ast h)$ live in the
vector space

\begin{equation}
H_{f, g, h} := \bigoplus_n H_{n;f, g, h} := \bigoplus_n
\left\langle \widetilde{X}^r f \widetilde{X}^s g \widetilde{X}^t
h\, \, \hbar^{r+s+t},  {r+s+t=n}\right\rangle.
\end{equation}

Generically, $H_{n;f, g, h}$ is a vector space of dimension
$\frac{1}{2} (n+1)(n+2)$. So it is natural to check the
identification of the coefficients with respect to the canonical
base $\widetilde{X}^r f \widetilde{X}^s g \widetilde{X}^t h\, \,
\hbar^{r+s+t}$ (${r+s+t=n}$). The problem is that in this case,
for $H_{n;f, g, h}$, we will have
$\displaystyle\sum_{r=0}^n\sum_{s=0}^{n-r}\sum_{t=0}^{n-r-s} 1=
\frac{1}{2} (n+1)(n+2)$ equations, which is not very practical.

In order to reduce the number of equations to verify, we will try
to determine a subspace in which live  $(f\ast g)\ast h$ and
$f\ast (g\ast h)$. In fact, we have already seen that when
 $f$ and $g$ are both holomorphic,  $f\ast g$ is a series which can be written
as a sum(with coefficients) of the $\hbar^n \sum (-1)^r {n\choose
r} \widetilde{\partial}^r f \widetilde{\partial}^{n-r} g$'s,  and
the latter form a basis of the kernel of the operator
$\hbar^{-1}\Delta E_-$,  we have then

\begin{lemma}
For three holomorphic functions $f, g, h\in \widetilde{\mathcal
M}$, the kernel of the operator $\hbar^{-1} E_-: H_{f, g,
h}\rightarrow H_{f, g, h}$ generated by the vectors ($0\leq p\leq
n$)

\begin{eqnarray*}
\xi_{n, p}& = & \hbar^n \sum_{s=0}^p (-1)^s {p\choose s}
\frac{\widetilde{X}^s}{(2k+2l+2n)_s} \left(\sum_{r=0}^{n-p}
{n-p\choose r} \widetilde{\partial}^{n-p-r} f
\widetilde{\partial}^{r} g\right) \widetilde{\partial}^{p-s} h.
\end{eqnarray*}

$(f\ast g)\ast h$ and $f\ast (g\ast h)$ belong to this kernel.

\end{lemma}

\noindent {\bf Proof.} We know that $H_{n;f, g, h}$ is a vector
space of dimension $\displaystyle \frac{1}{2} (n+1)(n+2)$. We
establish first the fact that the map ${\mathcal E}$ is
surjective: for every vector $\hbar^{n-1}\widetilde{\partial}^r f
\widetilde{\partial}^{s} g \widetilde{\partial}^t h$ with
$r+s+t=n-1$, we have

\begin{eqnarray*}
&& \hbar^{-1} E_-\Bigg(\hbar^n \sum_{i=0}^{n-1-r} \frac{(-1)^i
i!}{\prod_{u=0}^i (r+1+u)}\cr & & \hspace{1.2cm}
\left[\sum_{j=0}^i {s\choose i-j}{t\choose
j}\widetilde{\partial}^{r+1+i} f  \widetilde{\partial}^{s-i+j} g
\widetilde{\partial}^{t-j} h \right]\Bigg)\cr &=&
\frac{1}{2\pi}\hbar^{n-1}\sum_{i=0}^{n-1-r} \frac{(-1)^i i!
}{\prod_{u=0}^i (r+1+u)}\Bigg[\sum_{j=0}^i {s\choose i-j}{t\choose
j}\cr & & \hspace{1cm} \Big((r+1+i)\widetilde{\partial}^{r+i} f
 \widetilde{\partial}^{s-i+j} g \widetilde{\partial}^{t-j}
h\cr & & \hspace{2cm}+({s-i+j})\widetilde{\partial}^{r+1+i} f
\widetilde{\partial}^{s-i+j-1} g \widetilde{\partial}^{t-j} h\cr &
&\hspace{3cm}+({t-j})\widetilde{\partial}^{r+1+i} f
\widetilde{\partial}^{s-i+j} g \widetilde{\partial}^{t-j-1} h\Big)
\Bigg]\cr &=& \frac{1}{2\pi}\hbar^{n-1}\sum_{i,
j}\Bigg[\frac{(-1)^i i! }{\prod_{u=0}^{i-1} (r+1+u)}{s\choose
i-j}{t\choose j} \cr & & \hspace{3cm}+
\frac{(-1)^{i-1}(i-1)!}{\prod_{u=0}^{i-1} (r+1+u)}{s\choose
i-j}{t\choose j}i \Bigg]\cr & & \hspace{7cm}
\widetilde{\partial}^{r+i} f \widetilde{\partial}^{s-i+j} g
\widetilde{\partial}^{t-j} h\cr &=& \frac{1}{2\pi}\hbar^{n-1}
\widetilde{\partial}^{r} f \widetilde{\partial}^{s} g
\widetilde{\partial}^{t} h.
\end{eqnarray*}

The dimension at degree $n-1$ is $\displaystyle
\frac{1}{2} n(n+1)$,  this implies that the dimension of the kernel
at degree $n$ is $n+1$.

The vectors $\xi_{n, p}$ are in the kernel of $\hbar^{-1} E_-$: we
verify first that for two functions $f$ and $g$ in the kernel of $E_-$,  we have

\begin{eqnarray*}
E_- \widetilde{X}(f g) &=& 4\deg(f g) f g.
\end{eqnarray*}

So by simple induction, we can get

\begin{eqnarray}
 & & E_-\frac{\widetilde{X}^s}{(2k+2l+2n)_s}
\left(\sum_{r=0}^{n-p} {n-p\choose r} \widetilde{\partial}^{n-p-r}
f \widetilde{\partial}^{r} g\right) \cr &=&
\frac{\widetilde{X}^{s-1}}{(2k+2l+2n)_{(s-1)}}
\left(\sum_{r=0}^{n-p} {n-p\choose r} \widetilde{\partial}^{n-p-r}
f \widetilde{\partial}^{r} g\right),
\end{eqnarray}

\noindent which implies

\begin{eqnarray*}
\hbar^{-1} E_- \xi_{n, p} &=& \frac{1}{2\pi}\hbar^{n-1} \Bigg[
\sum_{s=0}^p (-1)^s {p\choose s} s \widetilde{\partial}^{s-1}
\left(\sum_{r=0}^{n-p} {n-p\choose r} \widetilde{\partial}^{n-p-r}
f
 \widetilde{\partial}^{r} g\right)
\widetilde{\partial}^{p-s} h\cr & & + \sum_{s=0}^p (-1)^s
{p\choose s} E_-\widetilde{\partial}^s \left(\sum_{r=0}^{n-p}
{n-p\choose r} \widetilde{\partial}^{n-p-r} f
\widetilde{\partial}^{r} g\right) (p-s)
\widetilde{\partial}^{p-s-1} h\Bigg]\cr &=&0.
\end{eqnarray*}

Moreover, we can project $\xi_{n, p}$ on the component whose
second factor is  $g$ and we get

\begin{equation}
\widetilde{\partial}^{n-p} f  g  \widetilde{\partial}^{p} h.
\end{equation}

\noindent These functions are generically linearly independent. This proves
that the $(n+1)$ $\xi_{n, p}$'s constitute a basis of the kernel of
$\hbar^{-1} E_-$ at degree $n$. $\Box$

\smallskip

In general, for all element $f\in\widetilde{\mathcal M}$, we can
define, in the vector space generated by the basis
 $\{ \varphi_n =
\frac{1}{(\deg f)_n} \widetilde{X}^n f,  \, \, n\in{\mathbb N}\}$,
an operator $\widetilde{\partial}$ by the formulae
$\widetilde{\partial} \varphi_n=\varphi_{n+1}$,  then
(\ref{tilpartilx}) is still valid. We can then define an operator
$\hbar^{-1} E_- : H_{f, g, h}\rightarrow H_{f, g, h}$ by the
following formula:

\begin{eqnarray}
\hspace{-0.8cm}\hbar^{-1} E_- ( \widetilde{\partial}^r f
\widetilde{\partial}^s g \widetilde{\partial}^t h\, \,
\hbar^{r+s+t}) & = & (r\widetilde{\partial}^{r-1} f
\widetilde{\partial}^s g \widetilde{\partial}^t h+s
\widetilde{\partial}^r f \widetilde{\partial}^{s-1} g
\widetilde{\partial}^t h\cr\hspace{-0.8cm} & &\hspace{1.5cm}
+t\widetilde{\partial}^r f \widetilde{\partial}^s g
\widetilde{\partial}^{t-1} h) \hbar^{r+s+t-1},
\end{eqnarray}

Then the above argument works without any modification.

\medskip

So it is sufficient now to identify the coefficients of $\hbar^n
\widetilde{\partial}^p f  g  \widetilde{\partial}^{n-p} h$ to
obtain the associativity. In $(f\ast g)\ast h$, it is the sum of
the terms (for $n-r\geq p$)

\[ \frac{(-1)^ r A_{r}(2k, 2l)}{(2k)_r}{n-r \choose p}\frac{(-1)^{p-r}A_{n-r}(2k+2l+2r,  2m)}{(2k+2l+2r)_{n-p} (2m)_p }. \]

For $f\ast (g\ast h)$, it is the sum of the terms (for $s\leq p$)

\[\frac{(-1)^p A_{n-s}(2k,  2l+2m+2s)}{(2k)_{n-p}(2l+2m+2s)_p} {n-s \choose n-p}\frac{A_{s}(2l, 2m)}{(2m)_s }\]

So finally what we should verify is the following identities, for
$p=0, 1, \dots,  n$:

\begin{eqnarray}\label{ident}
&& \sum_{r=0} {n-r \choose p} \frac{A_{n-r}(2k+2l+2r,  2m)
A_{r}(2k, 2l) }{(2k+2l+2r)_{n-p-r} (2m)_p  (2k)_r  }\cr & =&
\sum_{s=0} {n-s \choose n-p} \frac{A_{n-s}(2k,  2l+2m+2s)
A_{s}(2l, 2m) }{(2k)_{n-p}(2l+2m+2s)_{p-s} (2m)_s}.
\end{eqnarray}

We first look at the simplest case, the identification of the
coefficient of $\hbar$. We need to verify

\begin{eqnarray*} & & A_1(2k+2l,
2m)\left(\displaystyle\frac{1}{2k+2l}(f_{2k+2}g_{2l}h_{2m}+f_{2k}g_{2l+2}h_{2m})-f_{2k}g_{2l}\displaystyle\frac{1}{2m}h_{2m+2}\right)\cr
& & +  A_1(2k,  2l)\left(\displaystyle\frac{1}{2k}
f_{2k+2}g_{2l}h_{2m} -
f_{2k}\displaystyle\frac{1}{2l}g_{2l+2}h_{2m}\right)\cr &=
&A_1(2k, 2l+
2m)\left(\displaystyle\frac{1}{2k}f_{2k+2}g_{2l}h_{2m}-\displaystyle\frac{1}{2l+2m}(f_{2k}g_{2l+2}h_{2m}+f_{2k}g_{2l}h_{2m+2})\right)\cr
& &+ A_1(2l,  2m)\left(
f_{2k}\displaystyle\frac{1}{2l}g_{2l+2}h_{2m} -
f_{2k}g_{2l}\displaystyle\frac{1}{2m}h_{2m+2}\right) .
\end{eqnarray*}

\noindent In other words,

\begin{eqnarray*}
\frac{1}{2k+2l} A_1(2k+2l,  2m)+\frac{1}{2k}A_1(2k,  2l)& = &
\frac{1}{2k}A_1(2k, 2l+ 2m), \cr \frac{1}{2k+2l}A_1(2k+2l,
2m)-\frac{1}{2l}A_1(2k,  2l)& =&\frac{1}{2l}A_1(2l,
2m)-\frac{1}{2l+2m}A_1(2k, 2l+ 2m), \cr -\frac{1}{2m}A_1(2k+2l,
2m)& = & -\frac{1}{2l+2m}A_1(2k, 2l+ 2m)-\frac{1}{2m}A_1(2l,  2m).
\end{eqnarray*}

It is obvious that $A_1(2k, 2l)=2k\cdot 2l$ verify these
equations.

\medskip

Then we pass to the next step,  the identification of the
coefficients of $\hbar^2$:

\begin{eqnarray}
\hspace{-1.2cm}\frac{A_2(2k+2l, 2m)}{2m(2m+1)}  & = &
\frac{A_2(2k, 2l+2m)}{(2l+2m)(2l+2m+1)}\cr\hspace{-1.2cm} & &
\hspace{0.5cm} +4kl+\frac{A_2(2l, 2m)}{2m(2m+1)},
\cr\hspace{-1.2cm} \frac{A_2(2k+2l, 2m)}{(2k+2l)2m}+ (2k+2l+2) 2l
&　=　&　\frac{A_2(2k, 2l+2m)}{2k(2l+2m)}+(2l+2m+2)2l,
\cr\hspace{-1.2cm} \frac{A_2(2k+2l, 2m)}{(2k+2l)(2k+2l+1)} + 4lm+
\frac{A_2(2k, 2l)}{2k(2k+1)} & = & \frac{A_2(2k,
2l+2m)}{2k(2k+1)}.
\end{eqnarray}

This system has a special solution:

\begin{equation}
A_2(2k, 2l)=\frac{1}{2}2k(2k+1)2l(2l+1),
\end{equation}

\noindent so we need to solve the homogeneous system:

\begin{eqnarray}
\hspace{-1.2cm}\frac{A_2(2k+2l, 2m)}{2m(2m+1)}  & = &
\frac{A_2(2k, 2l+2m)}{(2l+2m)(2l+2m+1)}+\frac{A_2(2l,
2m)}{2m(2m+1)}, \cr\hspace{-1.2cm} \frac{A_2(2k+2l,
2m)}{(2k+2l)2m} &　=　&　\frac{A_2(2k, 2l+2m)}{2k(2l+2m)},
\cr\hspace{-1.2cm} \frac{A_2(2k+2l, 2m)}{(2k+2l)(2k+2l+1)} +
\frac{A_2(2k, 2l)}{2k(2k+1)} & = & \frac{A_2(2k,
2l+2m)}{2k(2k+1)}.
\end{eqnarray}

We note by $\widetilde{A}(2k, 2l)$ the function which
$\displaystyle\frac{2k+2l+1}{4kl}A(2k, 2l)$,  the equations that
$\widetilde{A}(2k, 2l)$ satisfy are:

\begin{eqnarray}
&  & \widetilde{A}_2(2k+2l,
2m)\left(\frac{1}{2m+1}-\frac{1}{2k+2l+2m+1}\right)\cr & = &
\widetilde{A}_2(2k,
2l+2m)\left(\frac{1}{2l+2m+1}-\frac{1}{2k+2l+2m+1}\right)\cr & &
\hspace{3cm} +\widetilde{A}_2(2l,
2m)\left(\frac{1}{2m+1}-\frac{1}{2l+2m+1}\right), \cr & &
\widetilde{A}_2(2k+2l, 2m) =　\widetilde{A}_2(2k, 2l+2m), \cr & &
\widetilde{A}_2(2k+2l, 2m)\left(\frac{1}{2k+2l+1}
-\frac{1}{2k+2l+2m+1}\right)\cr & & \hspace{3cm}
 +
\widetilde{A}_2(2k,
2l)\left(\frac{1}{2k+1}-\frac{1}{2k+2l+1}\right) \cr & = &
\widetilde{A}_2(2k,
2l+2m)\left(\frac{1}{2k+1}-\frac{1}{2k+2l+2m+1}\right).
\end{eqnarray}

The first two equations indicate $\widetilde{A}_2(2l,
2m)=\widetilde{A}_2(2k+2l, 2m)$ for all $(2k, 2l, 2m)$,  and by
using once more the second equation, we get $\widetilde{A}_2(2l,
2m)=\widetilde{A}_2(2k+2l, 2m)= \widetilde{A}_2(2k, 2l+2m)$, i.e.,
$\widetilde{A}$ is a constant function. We then conclude that in
our situation the degree of freedom is one,  i.e.,  in the general
formula of $A_2$ we can introduce a parameter $c$ :

\begin{equation}\label{II.49}
A_2(2k, 2l)=\frac{1}{2}2k(2k+1)2l(2l+1)+c \frac{2k 2l}{2k+2l+1}.
\end{equation}

Now we study some properties of a sequence $A_n$ which defines an
associative product. We assume their existence (the examples of
Cohen-Manin-Zagier provide some) and we have

\begin{lemma}
 Assume their existence,  the $A_n$'s ($n\geq 3$) are
determined by $A_0,  A_1, \dots,  A_{n-1}$ and the associativity.
\end{lemma}

\noindent{\bf Proof.} Our aim is to determine the value of $A_n(2x, 2y)$
 for every pair $(x, y)\in {\mathbb
N}^2\setminus \{(0, 0)\}$ ($(0, 0)$ is not included because in
this case, for all $n\geq 1$,  $[f, g]_n=0$). The idea is very
simple, in order to do the identification of the coefficients of
$\hbar^n$,  we have $n+1$ equations,  indexed by $p$,  by
considering $2k, 2l, 2m$ as constants and assume that $A_i(i<n)$
are already known.

If $l > 0$,  there is, in these equations,  (at most) four
unknowns: $A_n(2k, 2l)$, $A_n(2l, 2m)$, $A_n(2k+2l, 2m)$,
$A_n(2k, 2l+2m)$. The first two appear only once each:
$p=0$ for $A_n(2k, 2l)$,  and $p=n$ for $A_n(2l, 2m)$.
When $n\geq 3$,  we take the two equations with $p=1$ and
$2$. The determinant of the linear equation system with $A_n(2k+2l, 2m)$ and $A_n(2k, 2l+2m)$
as unknown is

\begin{eqnarray}
& &  \displaystyle \det \left( \begin{array}{rl}
     \displaystyle {n\choose 1} \frac{1}{(2k+2l)_{n-1}(2m)_1}  &  \displaystyle {n\choose n-1} \frac{1}{(2k)_{n-1}(2l+2m)_1}\\
         & \\
      \displaystyle {n\choose 2} \frac{1}{(2k+2l)_{n-2}(2m)_2} &   \displaystyle {n\choose n-2} \frac{1}{(2k)_{n-2}(2l+2m)_2}
      \end{array}\right)\cr
& & \cr &   =  & {n\choose n-1}{n\choose n-2}
      \frac{1}{(2k+2l)_{n-2}(2m)_1(2k)_{n-2}(2l+2m)_1}\cr
& & \hspace{0.5cm}
\left(\frac{1}{(2k+2l+n-2)(2l+2m+1)}-\frac{1}{(2m+1)(2k+n-2)}\right)\cr
&   =  & {n\choose n-1}{n\choose n-2}
      \frac{1}{(2k+2l)_{n-2}(2m)_1(2k)_{n-2}(2l+2m)_1}\cr
& & \hspace{0.5cm}
\frac{-(2l)^2-(2l)(2k+2m+n-1)}{(2k+2l+n-2)(2l+2m+1)(2m+1)(2k+n-2)}\neq
      0,
\end{eqnarray}

\noindent following the fact that $l>0, n>2$, and that $k, m$ are all positive integers.

We can therefore obtain the value of $A_n(2x, 2y)$ for a pair
 $(2x, 2y)$ which can be expressed as $(2k+2l, 2m)$ or $(2k,  2l+2m)$
for a certain $l>0$ without any ambiguity. The lemma is proven.$\Box$

\medskip

Next,  we have the following lemma by induction:

\begin{lemma}
We have $A_n(2k,  2l)=A_n(2l,  2k)$ and $A_n(2k,  0)=0$.
\end{lemma}

\noindent\bf Proof. \rm We have already obtained $A_n(2k,
2l)=A_n(2l,  2k)$ and $A_n(2k,  0)=0$ for $n=0, 1, 2$. Assume  now
that this is valid for $0, 1, \dots, n-1$. When we consider the
associativity identity for three functions
$f\in\widetilde{\mathcal M}_{2m}, g\in \widetilde{\mathcal M
}_{2l}, h\in \widetilde{\mathcal M }_{2k}$,  (\ref{ident})
becomes, for all fixed $n$ and $p$,

\begin{eqnarray*}
&& \sum_{r=0} {n-r \choose p} \frac{A_{n-r}(2m+2l+2r,  2k)
A_{r}(2m, 2l) }{(2m+2l+2r)_{n-p-r} (2k)_p  (2m)_r  }\cr & =&
\sum_{s=0} {n-s \choose n-p} \frac{A_{n-s}(2m,  2l+2k+2s)
A_{s}(2l, 2k) }{(2m)_{n-p}(2l+2k+2s)_{p-s} (2k)_s}.
\end{eqnarray*}

If we exchange the indices $r$ and $s$,  and replace
$p$ by $n-p$,  we obtain,

\begin{eqnarray*}
&& \sum_{s=0} {n-s \choose n-p} \frac{A_{n-s}(2m+2l+2s,  2k)
A_{s}(2m, 2l) }{(2m+2l+2s)_{p-s} (2k)_{n-p}  (2m)_s  }\cr & =&
\sum_{r=0} {n-r \choose p} \frac{A_{n-r}(2m,  2l+2k+2r) A_{r}(2l,
2k) }{(2m)_{p}(2l+2k+2r)_{n-p-r} (2k)_r}.
\end{eqnarray*}

For  $0<p<n$,  the only different with respect to (\ref{ident}),
by using induction hypothesis, is that we've replaced $A_n(2k,
2l+2m)$(resp. $A_n(2k+2l, 2m)$) by $A_n(2l+2m, 2k)$(resp. $A_n(2m,
2k+2l)$). This implies that $A_n(2l+2m, 2k)$ and $A_n(2m, 2k+2l)$
satisfy the same linear equation system as $A_n(2k, 2l+2m)$ and
$A_n(2k+2l, 2m)$,  the previous lemma gives $A_n(2x,  2y)=A_n(2y,
2x)$.

When we take $l=0, k, m\neq 0$ in (\ref{ident}),
The identity $p=0$ is simplified as

\begin{equation*}
\sum_{r=0} \frac{A_{n-r}(2k+2r,  2m) A_{r}(2k, 0)
}{(2k+2r)_{n-p-r} (2m)_p  (2k)_r} = \frac{A_{n}(2k,
2m+2s)}{(2k)_{n}},
\end{equation*}

\noindent i.e.

\begin{equation*}
\frac{A_{n}(2k,  2m) }{(2k)_{n}}+\frac{A_{n}(2k, 0) }{ (2k)_n} =
\frac{A_{n}(2k,  2m)}{(2k)_{n}},
\end{equation*}

\noindent then we have $A_{n}(2k, 0)$,  the lemma is established.$\Box$

\medskip

When we write  $A_n$ as a polynomial of $2k, 2l$ and $c$, then
because that $A_0,  A_1$ are both of degree 0 in $c$, we conclude
by the above argument that

\begin{lemma}\label{polyc}
$A_n$ is a polynomial of degree
$\left[\displaystyle\frac{n}{2}\right]$ in $c$.
\end{lemma}

\smallskip

In \cite{CMZ}, they use only the modularity to construct the
invariant formal pseudodifferential operators.

%Ainsi un produit $\ast$
%associatif au sens faible peut \^etre d\'efini sur
%$\widetilde{\mathcal M}[[\hbar]]$ au lieu de ${\mathcal M}[[\hbar]]$.
%Ce qu'ils ont obtenu est une famille \`a un param\`etre de
%produits formels(cf. (\ref{mukappa})). D'apr\`es la proposition
%ci-dessus,  les m\^emes coefficients d\'efinissent un produit
%$\star$ associatif au sens fort. Mais le Lemme \ref{polyc} nous
%dit qu'il y a exactement un param\`etre \`a utiliser pour obtenir
%tous les produits formels $\star$ associatifs au sens fort. On
%conclut donc

\begin{theorem}
Cohen-Manin-Zagier have in fact found all associative formal
products of the form (\ref{prodformfaib}).
\end{theorem}

\begin{remark}
We underline here two facts:
\begin{enumerate}
\item numerically,  the parameter $c$ introduced in
(\ref{II.49}) equals to $-3+4\kappa-\kappa^2$ for the $\kappa$ in
(\ref{mukappa});
\item when we consider the restriction to classical modular forms, for every degree
the space ${\mathcal M}_{2k}$ is of finite dimension. Our argument
above does not work any more, so it is not excluded that other
formal products defined using Rankin-Cohen brackets exist at this
level.
\end{enumerate}
\end{remark}

\medskip

We give a proposition which shows that the multiplication
structure defined by the Eholzer product (or Rankin-Cohen product for Connes-Moscovici)
is somewhat ``finer'' than that defined by the usual product,
in fact, we have

\begin{proposition}\label{38}
Let $\Gamma$ be a congruence subgroup of $SL_2({\mathbb Z})$ such
that ${\mathcal M}(\Gamma)$ admits the unique factorization
property (for example $SL_2({\mathbb Z})$ itself), let $F_1$,
$F_2$, $G_1$, $G_2\in {\mathcal M}(\Gamma)$ such that

\begin{equation}
RC(F_1, G_1)=RC(F_2, G_2),
\end{equation}

\noindent as formal series in ${\mathcal M}(\Gamma)[[\hbar]]$,
then there exists a constant $C$ such that

\begin{equation}
F_1=C F_2,  G_2= C G_1.
\end{equation}
\end{proposition}

We prove first that

\begin{lemma}\label{fine}
Let $f\in {\mathcal M}_{2k},  g\in {\mathcal M}_{2l}, h\in
{\mathcal M}_{2m} $ be three modular forms such that $[fg,
h]_n=[f, gh]_n$ for all $n$,  then $l=0$,  i.e.,  $g$ is a constant function.
\end{lemma}

\noindent{\bf Proof of the lemma}. Our data satisfy
automatically $[fg, h]_0=[f, gh]_0$. As to the case
$n=1$,  we have

\[ (2k+2l)fg\frac{dh}{dz}-2m\frac{d(fg)}{dz}h=2k f \frac{d(gh)}{dz}-(2l+2m)\frac{df}{dz},  \]

\noindent which implies that

\[(2k+2m)\frac{1}{g}\frac{dg}{dz}=2l \left(\frac{1}{f}\frac{df}{dz}+\frac{1}{h}\frac{dh}{dz}\right), \]

\noindent in other words

\begin{equation}\label{n=1}
g^{k+m}= C^{ste} (fh)^l,
\end{equation}

\noindent for a non-zero constant. Now we write the Fourier
expansions of these three modular forms:

\begin{eqnarray}
f &=& \alpha_0+\alpha_1 q+\alpha_2 q^2+ \cdots,  \cr g &=&
\beta_0+\beta_1 q+\beta_2 q^2+\cdots,  \cr h &= &\gamma_0+\gamma_1
q+\gamma_2 q^2+\cdots.
\end{eqnarray}

\noindent where $q=\exp(2\pi i z)$. As  $\displaystyle
q\frac{d}{dq}=\frac{1}{2\pi i}\frac{\partial}{\partial z}$,  we have

\begin{eqnarray}
\left(\frac{1}{2\pi i}\frac{\partial}{\partial z}\right)^n f &=&
0^ n\alpha_0+1^n \alpha_1 q+2^n \alpha_2 q^2+ \cdots,  \cr
\left(\frac{1}{2\pi i}\frac{\partial}{\partial z}\right)^n g &=&
0^n \beta_0+1^n \beta_1 q+2^n \beta_2 q^2+\cdots,  \cr
\left(\frac{1}{2\pi i}\frac{\partial}{\partial z}\right)^n h &= &
0^n \gamma_0+1^n \gamma_1 q+2^n \gamma_2 q^2+\cdots.
\end{eqnarray}

This implies that in the calculation of
$[fg, h]_n$ and of $[f, gh]_n$($n\geq 1$),  there are only two terms (among the $n+1$ to sum up)
which contain the term of degree 1 in $q$: the first and the last
in the definition formula. We have then for all $n$,

\begin{eqnarray}\label{n}
\hspace{-1.2cm}& & {2k+2l+n-1\choose n} \alpha_0\beta_0 \gamma_1
+(-1)^n{2m+n-1\choose n} (\alpha_0\beta_1+\alpha_1\beta_0)\gamma_0
\cr\hspace{-1.2cm} &=& {2k+n-1\choose n}
\alpha_0(\beta_0\gamma_1+\beta_1\gamma_0)+(-1)^n {2l+2m+n-1\choose
n} \alpha_1\beta_0\gamma_0.
\end{eqnarray}

We have to distinguish several different cases:

\medskip

\noindent 1) $l=0$,  i.e.,  $g=\beta_0$. (\ref{n}) is
automatically valid, and it's exactly the claim of the lemma.\\

\noindent 2) $l > 0$,  there are two possibilities:\\

a) $\beta_0\neq 0$,  then following (\ref{n=1}) we have
    $\alpha_0\neq 0, \gamma_0\neq 0$(because that the constant term of $(fh)^{k+m}$ is non-zero).
By using the bilinearity of the brackets,  it's possible to assume $\alpha_0=\beta_0=\gamma_0=1$.
    Then (\ref{n}) becomes,  for all $n$,
    \begin{eqnarray}
& & {2k+2l+n-1\choose n} \gamma_1 +(-1)^n{2m+n-1\choose n}
(\beta_1+\alpha_1) \cr &=& {2k+n-1\choose n}
(\gamma_1+\beta_1)+(-1)^n {2l+2m+n-1\choose n} \alpha_1.
\end{eqnarray}
Without loss of generality,  we can assume that $m\geq
k$(otherwise we consider $[hg, f]_n=[h, gf]_n$),  now the
variables $\alpha_1, \beta_1, \gamma_1$ satisfy the equations
(for all $n$)

\begin{equation}
A_{n1}\alpha_1+A_{n2}\beta_1+A_{n3}\gamma_1=0
\end{equation}

\noindent where

\begin{eqnarray}
A_{n1} & = & (-1)^n (2m)_n- (-1)^n (2l+2m)_n, \cr A_{n2} &=&
(-1)^n (2m)_n-(2k)_n,  \cr A_{n3} &=& (2k+2l)_n-(2k)_n,
\end{eqnarray}

\noindent and especially,

\begin{eqnarray}
A_{11}& = &  2l, \cr A_{12} & = & -2k-2m,  \cr A_{13} & = & 2l,
\cr A_{21} & = & -(2m+2m+1)2l -(2l)^2=-(4m+2l+1)(2l), \cr A_{22} &
= & 2m(2m+1)-2k(2k+1)=(2m-2k)(2m+2k+1),  \cr A_{23} & = &
(4k+2l+1)(2l), \cr A_{31} & = & (2l)^3+3(2m+1) (2l)^2+
(3(2m)^2+6(2m)+2)(2l),  \cr A_{32} & = &
-2m(2m+1)(2m+2)-2k(2k+1)(2k+2), \cr A_{33} & = & (2l)^3+3(2k+1)
(2l)^2+ (3(2k)^2+6(2k)+2)(2l).
\end{eqnarray}

The determinant of this system of linear equations is then

\begin{eqnarray}
& & \det A_{1\leq i, j \leq 3}\cr & = & \det \left(
                                     \begin{array}{lll}
                                      2l & -2k-2m & 2l \\
                                       -(4m+2l+1)(2l) & (2m-2k)(2m+2k+1) & (4k+2l+1)(2l) \\
                                       A_{31} & A_{32} & A_{33}\end{array} \right)\cr
    & = &  (2l)^2 \Big\{- 6(k+l+m+1)(2m-2k)^2\cr
         & &  +(2k+2l+2m+1)\Big[-(2k+2m)(6k+6l+2m)-12k\Big](2m-2k)\cr
         & & -(4k+4l+4m+2)(2k+2m)(2k+2l)(4k+2l+3)\Big\}\cr
         &<&0.
\end{eqnarray}

\noindent We have taken as hypothesis $m\geq k$ and as the weights
of modular forms $k, l, m$ are all positive integers and $l\geq
1$,  the last inequality is obtained because all the three terms
to be sum up are nonnegative.

We can conclude that $\alpha_1=\beta_1=\gamma_1=0$. The same
argument can be applied when we compare the coefficients of
$q^2$,  and we obtain a system of linear equations
for $\alpha_2, \beta_2, \gamma_2$ with the same coefficient matrix,
so $\alpha_2=\beta_2=\gamma_2=0$,  so on and so forth. we get a contradiction.\\

\noindent b)    $\beta_0=0$,  the argument in (\ref{n=1})
gives us
   $\alpha_0=0$ or $\gamma_0=0$.

    So we can then assume that the first nonzero terms are
   $\alpha_{r}q^r$,  $\beta_{s}q^s$,  $\gamma_{t}q^t$ ( $r, s, t\geq 0$ ). We consider now
  the term of lowest degree in $q$,  say $r+s+t$, in the
   identity $[fg, h]_n=[f, gh]_n$,  we obtain, for all $n$,

   \begin{eqnarray}\label{cas2-2}
    & &\sum_{p=0}^n (-1)^p {n\choose p}(2k+2l+p)_{n-p} (2m+n-p)_p (r+s)^p t^{n-p} \alpha_r\beta_s
    \gamma_t\cr
    &=& \sum_{q=0}^n (-1)^q {n\choose q}(2k+q)_{n-1} (2l+2m+n-q)_q r^q (s+t)^{n-q} \alpha_r\beta_s
    \gamma_t.
   \end{eqnarray}
   As the $\alpha_{r}, \beta_{s}, \gamma_{t}$ are nonzero, dividing both sides by $\alpha_r\beta_s
    \gamma_t$ this
   becomes
   \begin{eqnarray}
    & &\sum_{p=0}^n (-1)^p {n\choose p}(2k+2l+p)_{n-p} (2m+n-p)_p (r+s)^p t^{n-p} \cr
    &=& \sum_{q=0}^n (-1)^q {n\choose q}(2k+q)_{n-q} (2l+2m+n-q)_q r^q (s+t)^{n-q}.
       \end{eqnarray}
For $n=1$,  we have,
\[(k+m)s=l(r+t).\]
   By taking into account this relation, we obtain,  by replacing $s$ by $\displaystyle\frac{l(r+t)}{k+m}$
(if $k$ and $m$ are all zero,  then
   according to (\ref{n=1}),  $l=0$, too, a  contradiction) for every $n\geq
   2$ a homogeneous equation of degree $n$ in $r, t$. For $n=2$,  this equation is

   \begin{eqnarray*}
0 &= & [(2m)_2-(2l+2m)_2]r^2
-2[(2k+2l+1)(2m+1)-(2k+1)(2l+2m+1)]rt\cr & &
+[(2k+2l)_2-(2k)_2]t^2+2[(2m)_2+(2k+1)(2l+2m+1)]rs\cr & &
-2[(2k+2l+1)(2m+1)+(2k)_2]st+ [(2m)_2-(2k)_2]s^2\cr &=&
\frac{2l}{(k+m)^2} \Big\{\left[(k+3m)(k+l+m)+(k+m)\right]r^2 \cr &
&\hspace{0.5cm}
+(2m-2k)(k+l+m)rt-\left[(3k+m)(k+l+m)+(k+m)\right]t^2\Big\} .
   \end{eqnarray*}

   We see first that  $r$ and $t$ are either all zero or all non-zero,  because that the coefficients of $r^2$ and $t^2$ are all strictly non zero.
    The case where $r=t=0$ is already treated above, we assume from now on $r, s, t>0$.
The last expression has a factor $r+t$,  i.e.

\begin{eqnarray}
0&=& \frac{2l}{(k+m)^2} \Big\{\left[(k+3m)(k+l+m)+(k+m)\right]r^2
\cr & &\hspace{0.5cm}
+(2m-2k)(k+l+m)rt-\left[(3k+m)(k+l+m)+(k+m)\right]t^2\Big\}\cr &=&
\frac{2l}{(k+m)^2}(r+t)\Big\{[(k+3m)(k+l+m)+(k+m)]r\cr & &
\hspace{3cm} -[(3k+m)(k+l+m)+(k+m)]t\Big\}.
\end{eqnarray}

    This implies that there exists a positive constant $\mu$ such that

   \begin{eqnarray}\label{n=2res}
   t & = & \mu [(k+3m)(k+l+m)+(k+m)],  \cr
   r & = & \mu [(3k+m)(k+l+m)+(k+m)],  \cr
   s & = & \mu l [4(k+l+m)+2].
   \end{eqnarray}

We calculate the equation for $n=3$,  the difference of two sides is,  by using (\ref{n=1}),

\begin{eqnarray}
& & \frac{1}{(k+m)^3}\Big\{(2k + 2l)(2k + 2l + 1)(2k + 2l + 2)t^3
(k + m)^3\cr & & \hspace{1cm}-
  3(2k + 2l + 1)(2k + 2l + 2)(2m + 2)[(k + m)r + l(r + t)]t^2(k + m)^2 \cr & & \hspace{1cm}+
  3(2k + 2l + 2)(2m + 1)(2m + 2)[(k + m)r + l(r + t)]^2t(k + m) \cr & & \hspace{1cm}-
  2m(2m + 1)(2m + 2)[(k + m)r + l(r + t)]^3\cr & & \hspace{1cm} -
  2k(2k + 1)(2k + 2)[l(r + t) + t(k + m)]^3 \cr & &\hspace{1cm} +
  3(2k + 1)(2k + 2)(2l + 2m + 2)[l(r + t) + t(k + m)]^2 r(k + m)\cr & & \hspace{1cm} -
  3(2k + 2)(2l + 2m + 1)(2l + 2m + 2)[l(r + t) +
        t(k + m)] r^2(k + m)^2 \cr & & \hspace{1cm}+ (2l + 2m)(2l + 2m + 1)(2l + 2m +
        2)r^3(k + m)^3\Big\}
\end{eqnarray}

We denote by $P_3$ the braced quantity,  as an integer coefficient
polynomial of $k, l, m, r, t$. Taking the values of $r$ and $t$ as
in (\ref{n=2res}), we obtain a polynomial in  $k, l, m$ whose
coefficients are all {\bf positive} (cf. Appendix \ref{P3} for the
explicit expressions),  which implies that it could not have
positive integer roots in $k, l, m$. So this possibility is
excluded. $\Box$

\medskip

\noindent \bf Proof of the Proposition \ref{38}. \rm We do first a simplification: Let

\[F_1=f_{1, 2k}+f_{1, 2k+2}+f_{1, 2k+4}+\cdots \, \,  ;\, \,  G_1=g_{1, 2l}+g_{1, 2l+2}+g_{1, 2l+4}+\cdots ;\]

\[F_2=f_{2, 2k'}+f_{2, 2k'+2}+f_{2, 2k'+4}+\cdots \, \,  ;\, \,  G_2=g_{2, 2l'}+g_{2, 2l'+2}+g_{2, 2l'+4}+\cdots ;\]

\noindent be the natural graduation of these modular forms. Then when we look at, for each degree in $\hbar$,
the term whose coefficient is a modular form of smallest weight, we find the terms $[f_{1, 2k}, g_{1, 2l}]_n\hbar^n$ and $[f_{1,
2k'}, g_{1, 2l'}]_n\hbar^n$. So we have

\[RC(f_{1, 2k}, g_{1, 2l})  = RC(f_{2, 2k'},  g_{2, 2l'}).\]

In other words,  $[f_{1, 2k}, g_{1, 2l}]_n=[f_{2, 2k'}, g_{2,
2l'}]_n$ for all $n$. Using the unique factorization hypothesis,
we can speak of the biggest common divisor of  $f_{1, 2k}$ and
$f_{2, 2k'}$(resp. $g_{1, 2l}$ and $g_{2, 2l'}$),  denoted by
$f_0$(resp. $g_0$). We see first that by adjusting constants it's
possible to have

\begin{eqnarray}
\displaystyle\frac{f_{1, 2k}}{f_0}=\frac{g_{2, 2l'}}{g_0}&=&A, \cr
\displaystyle\frac{f_{2, 2k'}}{f_0}=\frac{g_{1, 2l}}{g_0}&=&B.
\end{eqnarray}

We have then $[f_0 A,  B g_0]_i=[f_0 B,  A g_0]_i$ for all
$i$.  Moreover,  $A, B$ are prime between them as polynomials of the generators. We use then

\begin{lemma}
Let $f\in{\mathcal M}_{2k}, A\in{\mathcal M}_{2l},
B\in{\mathcal M}_{2m}, g\in{\mathcal M}_{2n}$ be four modular forms
such that

\[[fA, Bg]_i=[fB, Ag]_i, \]

\noindent for all $i$,  and that $A, B$ are prime between them as polynomials of the generators,
then either $A=1$,  or
$B=1$.
\end{lemma}

\noindent\bf Proof of the lemma. \rm For $i=1$,  by definition,

\begin{equation}
(k+l)fA \frac{d(Bg)}{dz}- \frac{d(fA)}{dz}(m+n) Bg = (k+m)fB
\frac{d(Ag)}{dz}- \frac{d(fB)}{dz}(l+n) Ag,
\end{equation}

\noindent i.e.

\begin{eqnarray}
& & (k+l)fA \left(\frac{dB}{dz}g+B\frac{dg}{dz}\right)-
\left(\frac{df}{dz}A+f\frac{dA}{dz}\right)(m+n) Bg\cr &=& (k+m)fB
\left(\frac{dA}{dz}g+A\frac{dg}{dz}\right)-
\left(f\frac{dB}{dz}+\frac{df}{dz}B\right)(l+n) Ag.
\end{eqnarray}

 We divide the terms by $fABg$ to obtain

\begin{equation}
(l-m)\left(\frac{1}{f}\frac{df}{dz}+\frac{1}{g}\frac{dg}{dz}\right)
= (k+2m+n)\frac{1}{A}\frac{dA}{dz}-(k+2l+n)
\frac{1}{B}\frac{dB}{dz},
\end{equation}

\noindent i.e.

\begin{equation}
(fg)^{l-m}=\frac{A^{k+2m+n}}{B^{k+2l+n}}.
\end{equation}

If $l\geq m$,  then the left hand side is a polynomial in the generators.
As $A, B$ are prime between them, we get
 $B=1$. If $l\leq m$ we get  $A=1$. The lemma is proved. $\Box$

\medskip

We can summarize the two lemmas above as follows:

\begin{lemma}\label{bloc}
For four non-zero modular forms,  $f_1\in {\mathcal
M}_{2l}$,  $g_1\in {\mathcal M}_{2k}$,  $f_2\in {\mathcal
M}_{2l'}$,  $g_2\in {\mathcal M}_{2k'}$,  if we have

\begin{equation}
[f_1, g_1]_n=[f_2, g_2]_n,
\end{equation}

\noindent for all $n$. Then $k=k'$, $l=l'$,  and there exists a non-zero constant  $C$ such that

\begin{equation}
f_1=C f_2,  \, \, \, \, \, \,  Cg_1=g_2.
\end{equation}
\end{lemma}

\noindent\bf Proof of the Proposition \ref{38}(continued) . \rm
Following the Lemma \ref{bloc},  we have $k=k', l=l'$ and the
existence of a constant $C$ such that

\[f_{1, 2k}=C f_{2, 2k}, \, \, \,  g_{2, 2l}=C g_{1, 2l}.\]

Then we pass to the next degree,  i.e. in the expansion of $RC(F_1, G_1)=RC(F_2, G_2)$, of every $\hbar^n$,
the term with second lowest weight coefficient (which is an element in ${\mathcal M}(\Gamma)$).
Besides $f_{1, 2k}=Cf_{2, 2k}$ and $Cg_{1, 2l}=g_{2, 2l}$,  the relevant terms in the expansion of $F_1, G_1, F_2, G_2$ are,  $f_{1, 2k+2}$,  $f_{2, 2k+2}$,  $g_{1, 2l+2}$, $g_{2,
2l+2}$. we have,  for all $n$,

\begin{equation*}
[f_{1, 2k}, g_{1, 2l+2}]_n + [f_{1, 2k+2}, g_{1, 2l}]_n = [f_{2,
2k}, g_{2, 2l+2}]_n + [f_{2, 2k+2}, g_{2, 2l}]_n.
\end{equation*}

\noindent i.e.,

\[ [f_{1, 2k}, C g_{1, 2l+2}-g_{2, 2l+2}]_n=[f_{1, 2k+2}-C f_{2, 2k+2},  g_{1, 2l}]_n,  \]

\noindent for all $n$ and the same constant $C$.  If $f_{1,
2k+2}-C f_{2, 2k+2}\in {\mathcal M}_{2k+2}$ and $C g_{1,
2l+2}-g_{2, 2l+2}\in {\mathcal M}_{2l+2}$ are non-zero, we can
apply once more the Lemma \ref{bloc} to get a contradiction. So
the only possibility that left is

\[f_{1, 2k+2}=C f_{2, 2k+2}\, \, \,  , \, \, \,  g_{2, 2l+2}=C g_{1, 2l+2}.\]

The rest is an induction procedure. If we have already $f_{1,
2k+2i}=C f_{2, 2k+2i}$, $ g_{2, 2l+2i}=C g_{1, 2l+2i}$ for $0\leq
i\leq p-1$,  then when we consider in $RC(F_1, G_1)=RC(F_2, G_2)$
the term who belongs to ${\mathcal M}_{2k+2l+2n+2p} \hbar^n$,  we
get an equality

\begin{equation}
\sum_i [f_{1, 2k+2i}, g_{1, 2l+2p-2i}]_n = \sum_i [f_{2, 2k+2i},
g_{2, 2l+2p-2i}]_n.
\end{equation}

Using the induction hypothesis,  it can be simplified to

\begin{equation}
[f_{1, 2k}, g_{1, 2l+2p}]_n+[f_{1, 2k+2p}, g_{1, 2l}]_n = [f_{2,
2k}, g_{2, 2l+2p}]_n+[f_{2, 2k+2p}, g_{2, 2l}]_n,
\end{equation}

\noindent or,  in an equivalent way,

\[ [f_{1, 2k}, C g_{1, 2l+2p}-g_{2, 2l+2p}]_n=[f_{1, 2k+2p}-C f_{2, 2k+2p},  g_{1, 2l}]_n, \]

\noindent for all $n$ and the same constant $C$.  If $f_{1,
2k+2p}-C f_{2, 2k+2p}\in {\mathcal M}_{2k+2p}$ and $C g_{1,
2l+2p}-g_{2, 2l+2p}\in {\mathcal M}_{2l+2p}$ are both non-zero,
the Lemma \ref{bloc} gives rise to a contradiction. So it's
possible for us to conclude that

\[f_{1, 2k+2p}=C f_{2, 2k+2p}, \, \,  g_{2, 2l+2p}=C g_{1, 2l+2p}.\]

The proposition is established.$\Box$

\appendix

\section{The value of $P_3$}\label{P3}
\small The following are results of calculus of Mathematica.

\begin{eqnarray*}
& & P_3(k, l, m, r, t)\cr &=&4l(r + t)( -3k^2 r^2 - 2 k^3 r^2 + 3
k l r^2 + 2 k l^2 r^2 -
    6 k m r^2 - 15 k^2 m r^2 - 3 k^3 m r^2 + 3 l m r^2  \cr &&\hspace{0.5cm}--
    9 k l m r^2 - 6 k^2 l m r^2
    3 k l^2 m r^2 - 3 m^2 r^2 - 24 k m^2 r^2 -
    15 k^2 m^2 r^2 - 9 l m^2 r^2\cr &&\hspace{0.5cm} - 24 k l m^2 r^2 -
    9 l^2 m^2 r^2 - 11 m^3 r^2 - 21 k m^3 r^2 -
    18 l m^3 r^2 - 9 m^4 r^2 + 12 k^2 r t + 17 k^3 r t \cr &&\hspace{0.5cm} +
    3 k^4 r t + 6 k l r t + 21 k^2 l r t + 6 k^3 l r t +
    4 k l^2 r t + 3 k^2 l^2 r t+ 24 k m r t +
    51 k^2 m r t + 24 k^3 m r t\cr &&\hspace{0.5cm} + 6 l m r t +
    42 k l m r t + 42 k^2 l m r t + 4 l^2 m r t +
    18 k l^2 m r t + 12 m^2 r t + 51 k m^2 r t +
    42 k^2 m^2 r t \cr &&\hspace{0.5cm}+ 21 l m^2 r t + 42 k l m^2 r t +
    3 l^2 m^2 r t + 17 m^3 r t + 24 k m^3 r t +
    6 l m^3 r t + 3 m^4 r t -3 k^2 t^2 \cr &&\hspace{0.5cm}- 11 k^3 t^2 -
    9 k^4 t^2 + 3 k l t^2 - 9 k^2 l t^2 - 18 k^3 l t^2 +
    2 k l^2 t^2 - 9 k^2 l^2 t^2 - 6 k m t^2 -
    24 k^2 m t^2 \cr &&\hspace{0.5cm}- 21 k^3 m t^2 + 3 l m t^2 -
    9 k l m t^2 - 24 k^2 l m t^2 + 2 l^2 m t^2 -
    3 k l^2 m t^2 - 3 m^2 t^2 - 15 k m^2 t^2\cr &&\hspace{0.5cm} -
    15 k^2 m^2 t^2 - 6 k l m^2 t^2- 2 m^3 t^2 -
    3 k m^3 t^2+ 2 l^2 m r^2).
\end{eqnarray*}

By taking the values $t=\mu[(k+3m)(k+l+m)+(k+m)],
r=\mu[(3k+m)(k+l+m)+(k+m)]$,  one gets

\begin{eqnarray*}
& & P_3(k, l, m, \mu[(3k+m)(k+l+m)+(k+m)],
\mu[(3k+m)(k+l+m)+(k+m)])\cr &=& \mu^3 (48 k^5 l + 320 k^6 l + 720
k^7 l + 672 k^8 l + 256 k^9 l +
    96 k^4 l^2 + 960 k^5 l^2 + 2976 k^6 l^2 + 3552 k^7 l^2 \cr & & \hspace{12pt}+
    1536 k^8 l^2 + 640 k^4 l^3 + 3792 k^5 l^3 +
    6624 k^6 l^3 + 3584 k^7 l^3 + 1536 k^4 l^4 +
    5280 k^5 l^4 \cr & & \hspace{12pt}+ 4096 k^6 l^4+ 1536 k^4 l^5 +
    2304 k^5 l^5 + 512 k^4 l^6 + 240 k^4 l m + 1920 k^5 l m +
    5232 k^6 l m \cr & & \hspace{12pt}+ 5760 k^7 l m + 2304 k^8 l m +
    384 k^3 l^2 m + 4800 k^4 l^2 m + 18240 k^5 l^2 m +
    26016 k^6 l^2 m \cr & & \hspace{12pt}+ 12288 k^7 l^2 m + 2560 k^3 l^3 m +
    19152 k^4 l^3 m + 40896 k^5 l^3 m + 25088 k^6 l^3 m +
    6144 k^3 l^4 m \cr & & \hspace{12pt}+ 26784 k^4 l^4 m + 24576 k^5 l^4 m +
    6144 k^3 l^5 m + 11520 k^4 l^5 m + 2048 k^3 l^6 m +
    480 k^3 l m^2 \cr & & \hspace{12pt}+ 4800 k^4 l m^2 + 16080 k^5 l m^2 +
    21120 k^6 l m^2 + 9216 k^7 l m^2 + 576 k^2 l^2 m^2 +
    9600 k^3 l^2 m^2 \cr & & \hspace{12pt}+ 46176 k^4 l^2 m^2 +
    80352 k^5 l^2 m^2 + 43008 k^6 l^2 m^2 +
    3840 k^2 l^3 m^2 + 38496 k^3 l^3 m^2 \cr & & \hspace{12pt}+
    103968 k^4 l^3 m^2 + 75264 k^5 l^3 m^2 +
    9216 k^2 l^4 m^2 + 53952 k^3 l^4 m^2 +
    61440 k^4 l^4 m^2 \cr & & \hspace{12pt}+ 9216 k^2 l^5 m^2 +
    23040 k^3 l^5 m^2 + 3072 k^2 l^6 m^2 + 480 k^2 l m^3 +
    6400 k^3 l m^3 + 27120 k^4 l m^3 \cr & & \hspace{12pt}+ 43392 k^5 l m^3 +
    21504 k^6 l m^3 + 384 k l^2 m^3 + 9600 k^2 l^2 m^3 +
    61824 k^3 l^2 m^3 + 135840 k^4 l^2 m^3 \cr & & \hspace{12pt}
    +86016 k^5 l^2 m^3+ 2560 k l^3 m^3 + 38496 k^2 l^3 m^3 +
    139392 k^3 l^3 m^3 + 125440 k^4 l^3 m^3 \cr & & \hspace{12pt}+
    6144 k l^4 m^3 + 53952 k^2 l^4 m^3 + 81920 k^3 l^4 m^3 +
    6144 k l^5 m^3 + 23040 k^2 l^5 m^3\cr & & \hspace{12pt} + 2048 k l^6 m^3 +
    240 k l m^4 + 4800 k^2 l m^4 + 27120 k^3 l m^4 +
    54720 k^4 l m^4 + 256 l m^9\cr & & \hspace{12pt}+ 32256 k^5 l m^4 + 96 l^2 m^4 +
    4800 k l^2 m^4 + 46176 k^2 l^2 m^4 +
    135840 k^3 l^2 m^4 \cr & & \hspace{12pt}+ 107520 k^4 l^2 m^4 + 640 l^3 m^4 +
    19152 k l^3 m^4 + 103968 k^2 l^3 m^4 +
    125440 k^3 l^3 m^4 \cr & & \hspace{12pt}+ 1536 l^4 m^4 + 26784 k l^4 m^4 +
    61440 k^2 l^4 m^4 + 1536 l^5 m^4 + 11520 k l^5 m^4 +
    512 l^6 m^4 \cr & & \hspace{12pt}+ 48 l m^5 + 1920 k l m^5 + 16080 k^2 l m^5 +
    43392 k^3 l m^5 + 32256 k^4 l m^5 + 960 l^2 m^5 \cr & & \hspace{12pt}+
    18240 k l^2 m^5 + 80352 k^2 l^2 m^5 +
    86016 k^3 l^2 m^5 + 3792 l^3 m^5 + 40896 k l^3 m^5 \cr & & \hspace{12pt}+
    75264 k^2 l^3 m^5 + 5280 l^4 m^5 + 24576 k l^4 m^5 +
    2304 l^5 m^5 + 320 l m^6 + 5232 k l m^6 \cr & & \hspace{12pt}+
    21120 k^2 l m^6 + 21504 k^3 l m^6 + 2976 l^2 m^6 +
    26016 k l^2 m^6 + 43008 k^2 l^2 m^6 \cr & & \hspace{12pt}+ 6624 l^3 m^6 +
    25088 k l^3 m^6 + 4096 l^4 m^6 + 720 l m^7 +
    5760 k l m^7 + 9216 k^2 l m^7 \cr & & \hspace{12pt}+ 3552 l^2 m^7 +
    12288 k l^2 m^7 + 3584 l^3 m^7 + 672 l m^8 +
    2304 k l m^8 + 1536 l^2 m^8 ).
\end{eqnarray*}

\medskip

Projet AO, Institut de Math\'ematiques de Jussieu ,  175 rue du
Chevaleret 75013 Paris, France\\

email:yao@math.jussieu.fr

\end{document}